\documentclass[11pt,a4paper]{elsart}
\journal{}

\usepackage[centertags]{amsmath}
\usepackage{amsfonts}
\usepackage{amssymb}
\usepackage{newlfont}
\usepackage{color}
\usepackage{tabularx}
\usepackage{graphicx}

\newtheorem{definition}{\sc \bf Definition}[section]
\newtheorem{theorem}{\sc \bf Theorem}[section]
\newtheorem{example}{\sc \bf Example}[section]
\newtheorem{lemma}{\sc \bf Lemma}[section]
\newtheorem{proposition}{\sc \bf Propostion}[section]

\begin{document}
\begin{frontmatter}
\title{A reproducing kernel Hilbert space approach in meshless collocation method}

\author[IKIU]{Babak Azarnavid\corauthref{cor}},
\author[IKIU]{Mahdi Emamjome},
\author[PUT]{Mohammad Nabati},
\author[IKIU]{Saeid Abbasbandy}

\corauth[cor]{Corresponding author: babakazarnavid@yahoo.com}

\address[IKIU]{Department of Mathematics, Imam Khomeini International University, Ghazvin, 34149-16818, Iran}
\address[PUT]{Department of Basic Sciences, Abadan Faculty of Petroleum
engineering, Petroleum University of Technology, Abadan, Iran}

\begin{abstract}
In this paper we combine the theory of reproducing kernel Hilbert
spaces with the field of collocation methods to solve boundary
value problems with special emphasis on reproducing property of
kernels. From the reproducing property of kernels we proposed a
new efficient algorithm to obtain the cardinal functions of a
reproducing kernel Hilbert space which can be apply conveniently
for multi-dimensional domains. The differentiation matrices are
constructed and also we drive pointwise error estimate of applying
them. In addition we prove the non-singularity of collocation
matrix. The proposed method is truly meshless and can be applied
conveniently and accurately for high order and also
multi-dimensional problems. Numerical results are presented for
the several problems such as second and fifth order two point
boundary value problems, one and two dimensional unsteady Burgers'
equations and a parabolic partial differential equation in three
dimensions. Also we compare the numerical results with those
reported in the literature to show the high accuracy and
efficiency of the proposed method.
\end{abstract}

\end{frontmatter}

\section{Introduction}

Meshless methods for numerical solution of boundary value problems
have recently become more and more popular and several meshless
methods were proposed by authors. Meshless collocation methods based
on radial basis functions have been applied successfully in many
kind of differential equations \cite{1,2,3,4,24,14}. The meshless
collocation radial basis function methods are constructed in a
symmetric and an unsymmetric form. In both cases, a given partial
differential equation and various boundary conditions, are
discretized by point evaluations of both sides in certain
collocation nodes. The non-symmetric collocation method used by
Kansa is relatively simple to implement, however it results in an
unsymmetric system of equations due to the main differential
equation and boundary conditions and forcing these equations to be
satisfied at the grid points and the resulting unsymmetric
collocation matrix can be singular in exceptional cases \cite{4}.
Fasshauer \cite{14} shows that many of the algorithms and strategies
used for solving partial differential equations with polynomial
pseudospectral methods can be adapted for the use with radial basis
functions. In the most radial basis function methods for solving
differential equations, the reproducing property of kernels was not
used as characteristic property in reconstruction scheme. In this
paper we describe how meshless collocation method based on
reproducing kernel Hilbert spaces can be used to solve boundary
value problems numerically with emphasis on reproducing property of
kernels and without radial property. The collocation method is
usually implemented in the physical space by seeking approximate
solution $u$ of differential equation in the form
\begin{equation}\label{e1}
u_{N}(x)=\sum_{k=1}^{N}u_{N}(x_{k})h_{k}(x),
\end{equation}
where the functions $h_{k}$ are cardinal with respect to selected
collocation point set $X=\{x_{1},...,x_{N}\}$, that is
$h_{k}(x_{j})=\delta_{k,j}$ for $1\leq k,j\leq N$. Here we
consider the collocation method that uses the cardinal functions
belong to a reproducing kernel Hilbert space $(RKHS)$. One of the
advantage of radial basis functions is in complex geometries and
non-uniform discretizations, which has been kept in meshless
collocation Reproducing kernel Hilbert space method and in
addition the boundary points will not appear in discretizations
and collocation matrices. Instead of imposing the boundary
condition in collocation matrix, we try to use the reproducing
kernels that satisfies the homogenized boundary conditions of the
equation and then we can show the collocation matrix is
nonsingular if we have invertible bounded linear operator in main
equation. A new and convenient algorithm is proposed based on the
reproducing property of kernels, to obtain the cardinal functions
that can be applied for multi-dimensional case as simple as one
dimensional case. The differentiation matrix will be calculate and
we drive the pointwise error of its application to differentiation
in physical space. The error estimate of the reconstruction scheme
is given as a percentage of the norm of the true solution, which
is the only unknown quantity here and the error bound is based on
the interpolation error of a known function in reproducing kernel
Hilbert space $H$ which can be determined exactly for a given set
of collocation point. To demonstrate the efficiency and
performance of our method computationally, we apply it on some
problems such as high order differential equations and
multi-dimensional nonlinear problems and compare the results with
reported results in literature.

The rest of this paper is organized as follows: In Section 2 we
describe the construction of cardinal functions for a subspace of
reproducing kernel Hilbert space spanned by some translated kernels.
Differentiation matrices and the error estimate of employing them is
discussed in Section 3. In Section 4 we illustrate how the $RKHS$
collocation method can be applied on differential equations by
employing two linear differential equations, one and two dimensional
nonlinear Burgers' equations and a three dimensional time dependent
boundary value problem as test problems. Finally a brief conclusion
is given in Section 5.

\section{Cardinal functions}
The discussion of this section is devoted to the obtaining of the
multivariate cardinal functions of finite dimensional subspace of a
reproducing kernel Hilbert space spanned by translated kernels.

\begin{theorem} \label{th01} \cite{9}
Suppose that $H$ is a reproducing kernel Hilbert space with
reproducing kernel $K:\Omega \times \Omega\rightarrow \mathbb{R}$.
Then $K$ is positive semi definite. Moreover, $K$ is positive
definite if and only if the point evaluation functionals are
linearly independent in $H^{*}$.
\end{theorem}

In fact the $RKHS$ and its reproducing kernel define each other
uniquely. For a determined reproducing kernel $K:\Omega \times
\Omega\rightarrow \mathbb{R}$ of Hilbert space $H$ and
$X=\{x_{1},...,x_{N}\}\subset \Omega$, where $\Omega\subset
\mathbb{R}^{d}$, we define the finite dimensional linear
reconstruction space
\begin{equation}\label{e4}
\mathcal{S}_{X}=span\{ K(.,x_{j}):1\leq j \leq N\}.
\end{equation}
To construct the native space of the kernel from its reproducing
kernel $K$, we first define, the reconstruction space
\begin{equation}\label{e5}
\mathcal{S}=\{ s\in\mathcal{S}_{X}: X\subset
\mathbb{R}^{d},|X|<\infty \}
\end{equation}
containing all potential interpolants of the form
\begin{equation}\label{e6}
s(x)=\sum_{j=1}^{N}c_{j}K(x,x_{j}),
\end{equation}
for some $c=(c_{1},...,c_{n})^{T} \in \mathbb{R}^{N}$ and
$X=\{x_{1},...,x_{N}\}\subset \mathbb{R}^{d}$. Suppose that we are
given a vector
$$f_{X}=(f(x_{1}),...,f(x_{N}))^{T}\in \mathbb{R}^{N}$$
of discrete function values, sampled from an unknown function
$f:\mathbb{R}^{d}\rightarrow\mathbb{R}$ at a finite point set
$X=\{x_{1},...,x_{N}\}\subset\mathbb{R}^{d}$. The given data
$f_{X}$ has a unique interpolant of the form Eq.~\eqref{e6}, if
the kernel $K$ is positive definite. In this paper we focuses on
positive definite kernels.  For any $X=\{x_{1},...,x_{N}\}\subset
\mathbb{R}^{d}$, from the positive definiteness of the kernel $K$,
the matrix $A_{X}=(K(x_{i},x_{j}))_{1\leq i,j\leq
N}\in\mathbb{R}^{N\times N}$ is positive definite. Let
$f=\sum_{j=1}^{N}c_{j}K(.,x_{j})$ and
$g=\sum_{j=1}^{N}d_{j}K(.,x_{j})$ belongs to $\mathcal{S}$, then
using the reproducing property of kernel
$(K(.,x_{i}),K(.,x_{j}))_{H}=K(x_{i},x_{j})$, we can define the
inner product and norm in $\mathcal{S}$ generated by the positive
definite matrix $A_{X}$ as follows
$$(f,g)_{K}:=(\sum_{j=1}^{N}c_{j}K(x,x_{j}),\sum_{j=1}^{N}d_{j}K(x,x_{j}))_{H}=
\sum_{i=1}^{N}\sum_{j=1}^{N}c_{i}d_{j}K(x_{i},x_{j})=c^{T}A_{X}d,$$
$$\|f\|_{K}:=\sqrt{(f,f)_{K}}=\sum_{i=1}^{N}\sum_{j=1}^{N}c_{i}c_{j}K(x_{i},x_{j})=c^{T}A_{X}c,$$
for any $X=\{x_{1},...,x_{N}\}\subset\mathbb{R}$ and all $c,d\in
\mathbb{R}^{N}$. By completion of the inner product space
$\mathcal{S}$ with respect to its norm $\|.\|_{K}$, we obtain the
$RKHS$ correspond to $K$, $H=\overline{\mathcal{S}}$. The
orthonormal system $\{\psi_{i}(x)\}_{i=1}^{N}$
 of $\mathcal{S}_{X}$ can be derived from Gram-Schmidt
 orthogonalization process of $\{K(.,x_{i})\}_{i=1}^{N}$,
\begin{equation}\label{e7}
\psi_{i}(x)=\sum_{k=1}^{i}\beta_{ik}K(.,x_{k})\hskip .5 cm
(\beta_{ii}>0,i=1,2,....).
\end{equation}

\begin{proposition} \label{pr1}
There is for any point set $X=\{x_{1},...,x_{N}\}\subset
\mathbb{R}^{d}$ a unique cardinal basis
$\{h_{1},...,h_{N}\}\subset\mathcal{S}_{X}$ satisfying
\begin{equation}\label{e08}
h_{j}(x_{k})=   \left\{%
\begin{array}{ll}
   1  \hskip .5 cm  j=k\\
    0 \hskip .5 cm  j\neq k\\
\end{array}, \hskip .7 cm (1\leq
j,k\leq N), \right.
\end{equation}
and they can be obtain as follows
\begin{equation}\label{e8}
h_{i}(x)=\sum_{k=1}^{N}\beta_{ki}\psi_{k}(x),\hskip .7 cm (1\leq i
\leq N), x\in \mathbb{R}^{d}.
\end{equation}
\end{proposition}

\textbf{Proof.} From the uniqueness of interpolant
$s:\mathbb{R}^{d}\rightarrow \mathbb{R}$ satisfying $s_{X}=f_{X}$,i.e.,
\begin{equation}\label{e9}
s(x_{j})=f(x_{j}),\hskip .5 cm 1\leq j\leq N,
\end{equation}
for any $h_{i},(1\leq i\leq N)$, we have a unique interpolant
$s_{i}\in \mathcal{S}_{X}$ that satisfies Eq.~\eqref{e08}, and can
be given as follows
\begin{equation*}\label{e09}
\begin{array}{ll}
  s_{i}(x)=\sum_{k=1}^{N}c_{ik}\psi_{k}(x)=\sum_{k=1}^{N}(s_{i},\psi_{k})_{H}\psi_{k}(x)= \\
 \sum_{k=1}^{N}(s_{i},\sum_{j=1}^{k}\beta_{kj}K(.,x_{j}))_{H}\psi_{k}(x)=\\
  \sum_{k=1}^{N}\sum_{j=1}^{k}\beta_{kj}(s_{i},K(.,x_{j}))_{H}\psi_{k}(x)=\\
    \sum_{k=1}^{N}\sum_{j=1}^{k}\beta_{kj}s_{i}(x_{j})\psi_{k}(x)=\\
  \sum_{k=1}^{N}\sum_{j=1}^{k}\beta_{kj}h_{i}(x_{j})\psi_{k}(x)=
  \sum_{k=1}^{N}\beta_{ki}\psi_{k}(x),\hskip .5 cm 1\leq i\leq N,  \\
\end{array}
\end{equation*}
where the functions $s_{i}\mid_{(1\leq i \leq N)}$ are the
involved cardinal functions $h_{i}\mid_{(1\leq i \leq N)} \in
\mathcal{S}_{X}$.

One of advantage of the above proposition is that the cardinal
functions can be easily construct for the multivariate case. In
\cite{25} the  cardinal functions of $\mathcal{S}_{X}$ are given
based on inversion of matrix $A_{X}$. We avoids inverting an ill
conditioned matrix and calculate the cardinal functions using the
reproducing property of kernels and a modified Gram Schmidt
orthogonalization algorithm, that cost $O(N^{2})$ flops versus
inversion of matrix $O(N^{3})$ flops.

\section{Differentiation matrices}
One of several ways to implement the collocation method is via so
called differentiation matrices. Suppose that $H$ is a reproducing
kernel Hilbert space with reproducing kernel $K:\Omega \times
\Omega\rightarrow \mathbb{R}$. When we just have values of a
function $u\in H$ at scattered point $X=\{x_{1},...,x_{N}\}$ in the
domain of $u$, and then the value $\mathcal{L}u(z)$ for a fixed
point $z\in \Omega$ must be approximated via
\begin{equation}\label{e010}
\mathcal{L}u(z)\simeq\sum_{k=1}^{N}\alpha_{k}u(x_{k}),
\end{equation}
where $\mathcal{L}$ is a linear bounded differential operator. Our
aim is to reconstruct operation of $\mathcal{L}$ on an unknown
function $u\in H$ from its values $u_{X}\in \mathbb{R}^{N}$. Define
the linear subspace of $H$,
$$\mathcal{S}_{X}=span\{ K(.,x_{j}):1\leq j \leq N\},$$
let $\mathcal{S}_{X}^{\perp}$ be the linear subspace of $H$,
$$\mathcal{S}_{X}^{\perp}=\{f\in H|f(x_{i})=0,i=1,...,N\}.$$
From the reproducing property of $K$, for any $f\in
\mathcal{S}_{X}^{\perp}$ we have
$$(f,\sum_{k=1}^{N}c_{k}K(.,x_{k}))_{H}=\sum_{k=1}^{N}c_{k}(f,K(.,x_{k}))_{H}=\sum_{k=1}^{N}c_{k}f(x_{k})=0,$$
then $\mathcal{S}_{X}^{\perp}$ is the orthogonal complement of
$\mathcal{S}_{X}$ and every $u\in H$ can be uniquely decomposed to
$u=u_{0}+u_{0}^{\perp}$, where $u_{0}\in\mathcal{S}_{X}$ and
$u_{0}^{\perp}\in\mathcal{S}_{X}^{\perp}$. For $1\leq k\leq N,$ we
have $u(x_{k})=u_{0}(x_{k})+u_{0}^{\perp}(x_{k})=u_{0}(x_{k})$ then
$u_{0}(x)=\sum_{k=1}^{N}c_{k}K(x,x_{k})=\sum_{k=1}^{N}u_{0}(x_{k})h_{k}(x)=\sum_{k=1}^{N}u(x_{k})h_{k}(x),$
is unique interpolant of $u\in S_{X}$. From Riesz representation
theorem, for any bounded linear operator $\mathcal{L}$ there exist
$g\in H$ such that $\mathcal{L}u(z)=(g,u)_{H}$ for a fixed point
$z\in \Omega$ and any $u\in H$, then
\begin{equation}\label{e011}
\begin{array}{ll}
\mathcal{L}u(z)=(g,u)_{H}\simeq (g,u_{0})_{H}=(g_{0}+g_{0}^{\perp},u_{0})_{H}=(g_{0},u_{0})_{H}=\\
(\sum_{k=1}^{N}\alpha_{k}K(.,x_{k}),u_{0})_{H}=\sum_{k=1}^{N}\alpha_{k}u_{0}(x_{k})=
\sum_{k=1}^{N}\alpha_{k}u(x_{k}).\\
\end{array}
\end{equation}
We can obtain the coefficients $\alpha_{i},(1\leq i\leq N)$ as
follows
\begin{equation*}\label{e012}
\begin{array}{ll}
\mathcal{L}h_{i}(z)=(g,h_{i})_{H}=(g_{0}+g_{0}^{\perp},h_{i})_{H}=(g_{0},h_{i})_{H}=(\sum_{k=1}^{N}
\alpha_{k}K(.,x_{k}),h_{i})_{H}=\\
\sum_{k=1}^{N}\alpha_{k}h_{i}(x_{k})=\alpha_{i}.\\
\end{array}
\end{equation*}
It is clear that Eq.~\eqref{e011} treats the functions
$u\in\mathcal{S}_{X}$ exactly.  If we evaluate Eq.~\eqref{e010} at
the grid point $x_{i},(i=1,...,N),$ for any $u\in\mathcal{S}_{X}$
then we get
$$\mathcal{L}u(x_{i})=\sum_{k=1}^{N}u(x_{k})\mathcal{L}h_{k}(x_{i}),\hskip .7 cm i=1,...,N$$
or in matrix vector notation $$\mathcal{L}\mathbf{u}=L\mathbf{u},$$
where
$\mathcal{L}\mathbf{u}=(\mathcal{L}u(x_{1}),...,\mathcal{L}u(x_{N}))^{T},
\mathbf{u}=(u(x_{1}),...,u(x_{N}))^{T}$ and the entries of the
differentiation matrix $L$ are given by
$l_{ij}=\mathcal{L}h_{j}(x_{i})$, which can be obtained as follows.
From \eqref{e8} we have
\begin{equation*}\label{e00_r1}
\begin{array}{ll}
l_{ij}&=\mathcal{L}h_{j}(x_{i})=\sum_{k=1}^{N}\beta_{kj}\mathcal{L}\psi_{k}(x_{i})\\
&=\left(\mathcal{L}\psi_{1}(x_{i}),\mathcal{L}\psi_{2}(x_{i}),...,
\mathcal{L}\psi_{N}(x_{i})\right).(\beta_{1j},\beta_{2j},...,\beta_{Nj})^{T}.\\
\end{array}
\end{equation*}
In practice the differentiation matrix can be obtained as follows
\begin{equation*}\label{e00_r2}
L=\Phi^{T}.B,
\end{equation*}
where $\Phi_{i,j}=\mathcal{L}\psi_{i}(x_{j})$ and $B$ is a lower
triangular matrix with $B_{i,j}=\beta_{ij}$. The differentiation
matrix $L$ can now be used to solve partial differential equations.
Sometimes only multiplication by $L$ is required e.g., for many time
dependent algorithm, and for other problem one needs to be able to
invert $L$.

\begin{lemma} \label{lem1}
If the linear operator $\mathcal{L}$ is invertible then the
corresponding matrix $L$ is nonsingular.
\end{lemma}

\textbf{Proof.} For selected collocation point set
$X=\{x_{1},...,x_{N}\}$, suppose that
$$c_{1}\mathcal{L}h_{1}(x)+...+c_{N}\mathcal{L}h_{N}(x)=0,$$ then from
the linearity of operator $\mathcal{L}$ we have
$$\mathcal{L}(c_{1}h_{1}(x)+...+c_{N}h_{N}(x))=0.$$
It follows that $$c_{1}h_{1}(x)+...+c_{N}h_{N}(x)\equiv0$$ from the
existence of $\mathcal{L}^{-1}$. Then from the linearly independence
of cardinal functions we conclude that the functions
$\mathcal{L}h_{j}(x)|_{1\leq j\leq N}$ are linearly independent and
the matrix $L$ has full rank.

\begin{theorem} \label{th03}
Let $H$ be a reproducing kernel Hilbert space with reproducing
kernel $K:\Omega\times\Omega\rightarrow \mathbb{R}$, and let
$\mathcal{L}$ be a linear differential operator. Assume that
$X=\{x_{1},...,x_{N}\}\subset \Omega$, where $\Omega \subset
\mathbb{R}^{d}$. Then for any $z\in \Omega$ and $u\in H$ we have
\begin{equation}\label{e12}
|\mathcal{L}u(z)-\sum_{k=1}^{N}u(x_{k})\mathcal{L}h_{k}(z)|\leq
\|u\|_{H}\|\varepsilon_{X}\|_{H},
\end{equation}
where
$\|\varepsilon_{X}\|_{H}=\|\mathcal{L}_{y}K(.,z)-\sum_{k=1}^{N}K(.,x_{k})\mathcal{L}h_{k}(z)\|_{H}$
is the norm of interpolation error of a known function
$\mathcal{L}_{y}K(.,z)$ and  $0\leq \|\varepsilon_{X}\|^{2}_{H} \leq
\mathcal{L}_{x}\mathcal{L}_{y}K(x,y)|_{x,y=z}$.
\end{theorem}

\textbf{Proof.} From Riesz representation theorem, for any bounded
linear operator $\mathcal{L}$ there exist $g\in H$ such that
$\mathcal{L}u(z)=(g,u)_{H}$ for any $u\in H$ and $z\in\Omega$, then
\begin{equation*}\label{e13}
\begin{array}{ll}
|\mathcal{L}u(z)-\sum_{k=1}^{N}u(x_{k})\mathcal{L}h_{k}(z)|=|(u,g)_{H}-(u,\sum_{k=1}^{N}K(.,x_{k})
\mathcal{L}h_{k}(z))_{H}|=\\
\|(u,g-\sum_{k=1}^{N}K(.,x_{k})\mathcal{L}h_{k}(z))_{H}\|\leq \|u\|_{H}\|g-\sum_{k=1}^{N}K(.,x_{k})\mathcal{L}h_{k}(z)\|_{H}.\\
\end{array}
\end{equation*}
Let
$\|\varepsilon_{X}\|_{H}=\|g-\sum_{k=1}^{N}K(.,x_{k})\mathcal{L}h_{k}(z)\|_{H}$.
From Eq.~\eqref{e011} it is easy to see that
$\mathcal{L}h_{k}(z)|_{1\leq k\leq N}$ are the interpolation
coefficients of $g=\mathcal{L}_{y}K(.,z)\in H$, and then
$\|\varepsilon_{X}\|_{H}$ is the interpolation error of known
function $\mathcal{L}_{y}K(.,z)$ which can be determined exactly.
Then
\begin{equation*}\label{e14}
\begin{array}{ll}
\|\varepsilon_{X}\|_{H}^{2}=(g-\sum_{k=1}^{N}K(.,x_{k})\mathcal{L}h_{k}(z),g-\sum_{k=1}^{N}K(.,x_{k})\mathcal{L}h_{k}(z))_{H}=\\
(g,g)_{H}-2\sum_{k=1}^{N}(g,K(.,x_{k}))_{H}\mathcal{L}h_{k}(z)+\sum_{k=1}^{N}\sum_{j=1}^{N}K(x_{k},x_{j})\mathcal{L}h_{k}(z)\mathcal{L}h_{j}(z)=\\
(g,g)_{H}-2\sum_{k=1}^{N}\sum_{j=1}^{N}K(x_{j},x_{k})\mathcal{L}h_{j}(z)\mathcal{L}h_{k}(z)+\sum_{k=1}^{N}\sum_{j=1}^{N}K(x_{k},x_{j})\mathcal{L}h_{k}(z)\mathcal{L}h_{j}(z)=\\
(g,g)_{H}-(\mathcal{L}h(z))^{T}A_{X}(\mathcal{L}h(z)),
\end{array}
\end{equation*}
where
$\mathcal{L}h(z)=(\mathcal{L}h_{1}(z),...,\mathcal{L}h_{N}(z))^{T}.$
So from positive definiteness of matrix $A_{X}$ and
$(g,g)_{H}=\mathcal{L}g(z)=\mathcal{L}(g,K(.,y))_{H}|_{y=z}=\mathcal{L}_{x}\mathcal{L}_{y}K(x,y)|_{x,y=z}$
we have
$$0\leq\|\varepsilon_{X}\|^{2}_{H}\leq\mathcal{L}_{x}\mathcal{L}_{y}K(x,y)|_{x,y=z},$$
where the subscript index shows the variable that the linear
operator acts on.
The above theorem gives us a pointwise error estimate of the
reconstruction scheme as a percentage of the norm of the true
solution, which is the only unknown quantity here and the error
bound is based on the interpolation error of a known function in
reproducing kernel Hilbert space $H$. Suppose we have a linear
differential equation of the form
\begin{equation}\label{e0011}
\mathcal{L}u=f,
\end{equation}
for some given $f\in V$ by ignoring boundary condition, where
$\mathcal{L}:U\rightarrow V$ is a linear data map that takes each
$u$ in some linear normed space $U$ into the data space $V$. We
use the following definition of well-posedness for above problem
given in \cite{26}.

\begin{definition} \label{de001}
An analytic problem as \eqref{e0011} is well-possed with respect to
a well-posedness norm $\|.\|_{W}$ on $U$ if there is a constant $C$
such that a well-posedness inequality
\begin{equation}\label{e0012}
\|u\|_{W}\leq C\|\mathcal{L}u\|_{V}, u\in U
\end{equation}
holds.
\end{definition}

This means that $\mathcal{L}^{-1}$ is continuous as a map
$\mathcal{L}(U)\rightarrow U$ in the norms $\|.\|_{V}$ and
$\|.\|_{W}$. We approximate $u^{*}$ the true solution of
\eqref{e0011} by $u^{*}_{N}$ defined as \eqref{e1} in finite
dimensional $U_{N}\subseteq U$, which can be calculate as a
finite-dimensional approximation problem from minimization the
residual norm $\|\mathcal{L}u_{N}-f\|_{V}$ over all $u_{N}\in
U_{N}$.

\begin{theorem} \label{th003}
Suppose the conditions of Theorem \ref{th03} hold. Assume a
well-posed analytic problem $\mathcal{L}u=f$ for some given $f\in
V$ in the sense of Definition \ref{de001} with some constant $C$,
well-posedness norm $\|.\|_{W}$ on some RKHS like $H$ and the data
space norm $\|.\|_{V}=max_{z\in \Omega}|.|$. If $u^{*}$ is the
true solution and $u^{*}_{N}\in\mathcal{S}_{X}\subseteq H $ is the
approximated solution, then we have
\begin{equation}\label{e0013}
\| u^{*}-u^{*}_{N}\|_{W}\ \leq C \|u\|_{H}\|\varepsilon_{X}\|_{H},
\end{equation}
where $\|\varepsilon_{X}\|_{H}$ is as defined in Theorem
\ref{th03}.
\end{theorem}

\textbf{Proof.} The true solution $u^{*}\in H$ can be uniquely
decomposed to $u^{*}=u_{0}+u_{0}^{\perp}$, where
$u_{0}\in\mathcal{S}_{X}$. From the well-posedness condition
\eqref{e0012} and \eqref{e12} we have
\begin{equation*}\label{e0014}
\begin{array}{ll}
\| u^{*}-u^{*}_{N}\|_{W}\leq C \| \mathcal{L}(u^{*}-u^{*}_{N})\|_{V}
\leq C \| \mathcal{L}(u^{*}-u_{0})\|_{V} \\=C\hskip .1 cm max_{z\in
\Omega}|\mathcal{L}(u^{*}-u_{0})(z)|=C\hskip .1 cm max_{z\in
\Omega}|\mathcal{L}u^{*}(z)-\mathcal{L}u_{0}(z)|\\
=C|\mathcal{L}u^{*}(z_{0})-\mathcal{L}u_{0}(z_{0})|
=C|\mathcal{L}u^{*}(z_{0})-\sum_{k=1}^{N}u_{0}(x_{k})\mathcal{L}h_{k}(z_{0})|\\
=C|\mathcal{L}u^{*}(z_{0})-\sum_{k=1}^{N}u^{*}(x_{k})\mathcal{L}h_{k}(z_{0})|\leq C \|u^{*}\|_{H}\|\varepsilon_{X}\|_{H}\\
\end{array}
\end{equation*}
where $z_{0}$ is a point in $\Omega$.

This theorem shows that if a problem in mathematical analysis is
well-posed, it have a discretization in numerical analysis that is
also well-posed \cite{26}. The Eq.~\eqref{e0013} describes the
worst–case error behavior of the approximated solution
$u^{*}_{N}$, and the error is given as a percentage of
$\|u^{*}\|_{H}$, which is the only unknown quantity in the error
bound.

\section{Implementation of the method}
We now begin with a discussion of some $RKHS$ and their reproducing
kernels for construction of differentiation matrices.

\begin{definition} \label{de01}
Let $H$ denote a Hilbert space of functions
$f:\Omega\subset\mathbb{R}^{d}\rightarrow \mathbb{R}$. Then a
function $K:\Omega \times \Omega\rightarrow \mathbb{R}$ is said to
be reproducing kernel for $H$, if and only if $K(.,x)\in H$ for all
$x\in \mathbb{R}^{d}$, and $(K(.,x),f)_{H}=f(x)$ for all $f\in H$
and all $x\in \mathbb{R}^{d}$.
\end{definition}

\begin{definition} \label{de1}
The inner product space $W_{2}^{m}[a,b]$ is defined as
$W_{2}^{m}[a,b]=\{ u(x)|u^{(m-1)}$ is absolutely continuous real
valued functions, $u^{(m)}\in L^{2}[a,b]\}$. The inner product in
$W_{2}^{m}[0,1]$ is given by
\begin{equation}\label{e2}
(u(.),v(.))_{W_{2}^{m}}=\sum_{i=0}^{m-1}u^{(i)}(a)v^{(i)}(a)+\int_{a}^{b}u^{(m)}(x)v^{(m)}(x)dx,
\end{equation}
and the norm $\|u\|_{W_{2}^{m}}$ is denoted by
$\|u\|_{W_{2}^{m}}=\sqrt{(u,u)_{W_{2}^{m}}}$,where $u,v\in
W_{2}^{m}[a,b]$.
\end{definition}

\begin{theorem} \label{th1} \cite{5}
The space $W_{2}^{m}[a,b]$ is a reproducing kernel space. That is,
for any $u(.)\in W_{2}^{m}[a,b]$ and each fixed $x\in[a,b]$, there
exists $R_{x}(.)\in W_{2}^{m}[a,b]$, such that
$(u(.),R_{x}(.))_{W_{2}^{m}}=u(x)$. The reproducing kernel
$R_{x}(.)$ can be denoted by
\begin{equation}\label{e3}
  R_{y}(x)=  \left\{%
\begin{array}{ll}
    \sum_{i=1}^{2m}c_{i}(y)x^{i-1} \hskip .1cm, & \hbox{$x\leq y$,} \\
    \sum_{i=1}^{2m}d_{i}(y)x^{i-1}\hskip .1cm, & \hbox{$x>y$}.\\
\end{array}%
\right.
\end{equation}
\end{theorem}

For more detail about reproducing kernel Hilbert spaces
$W_{2}^{m}[a,b]$ and the method of obtaining their reproducing
kernels $R_{x}(y)$, refer to \cite{5} and references there in.

\begin{proposition} \label{pr2}
The reproducing kernel of $W_{2}^{m}[a,b]$, for a bounded interval
$[a,b]$ is strictly positive definite.
\end{proposition}

\textbf{Proof.} Based on Theorems \ref{th1} and \ref{th01}, if we
can show that the point evaluation functionals are linearly
independent, the proof is complete. For pairwise distinct point
$x_{1},...,x_{N}\in [a,b]$ and $c\in \mathbb{R}^{N}$. Suppose that
$c_{1}\delta_{x_{1}}+...+c_{N}\delta_{x_{N}}=0$ then for any $f\in
W_{2}^{m}[a,b]$ we have
$$(c_{1}\delta_{x_{1}}+...+c_{N}\delta_{x_{N}})f=c_{1}f(x_{1})+...+c_{N}f(x_{N})=0.$$
Let $f(x)=x,x^{2},...,x^{N}\in W_{2}^{m}[a,b]$, then we have the
matrix system $Vc=0$ where $V$ is the $N\times N$ Vandermonde matrix
\begin{equation}\label{e15}
 \left(
   \begin{array}{ccc}
     x_{1}& \cdots & x_{N} \\
     x_{1}^{2}& \cdots & x_{N}^{N}\\
     \vdots& \ddots & \vdots \\
     x_{1}^{N}& \cdots & x_{N}^{N}\\
   \end{array}
 \right)c=0.
\end{equation}
From the nonsingularity of Vandermonde matrix we have $c=0$, so
the point evaluation functionals are linearly independent.

In the following we illustrate how the differentiation matrix can be
applied to solve differential equation. Several examples are
discussed in this section which provide samples of how simple
differentiation matrix can be applied to the boundary value problems
and high-dimensional time dependent problems. In order to solve
following problems, we may construct the closed subspaces
$\widehat{W}_{2}^{m}[a,b]$ of the reproducing kernel space
$W_{2}^{m}[a,b]$ by imposing several homogeneous boundary condition
on $\widehat{W}_{2}^{m}[a,b]$. To initially try the application of
differentiation matrix we consider linear two-point boundary value
problems of second and fifth order. Then one and two dimensional
Burgers' equations have been considered, as application of proposed
method on time dependent problems. Our final example is a three
dimensional problem, for explore the power of the method for solving
the multi-dimensional problems. To show the efficiency of the
present method for our problems in comparison with the exact
solution and other reported results, we report maximum absolute
error, the norm of relative errors and root mean squared error of
the solutions. In the proposed method, firstly, the nonhomogeneous
problem is reduced to a homogeneous one and then the functions
$\widehat{R}_{x_{j}}(x),j=1,...,N$ are used as the basis functions
to approximate the solution of the homogenized problem, where
$\widehat{R}_{y}(x)$ is the reproducing kernel of
$\widehat{W}_{2}^{m}[a,b]$, hence the approximate solution satisfies
the boundary conditions exactly. For the method of obtaining the
reproducing kernels $\widehat{R}_{y}(x)$, refer to \cite{5} and
references there in. Let
\begin{equation*}\label{bc_r1}
 Lu(\textbf{x})=f(\textbf{x}),\hspace{.5 cm} \textbf{x}\in \Omega\subset R^{d},\ Bu(\textbf{x})=g(\textbf{x}),\hspace{.5 cm} \textbf{x}\in \partial \Omega
\end{equation*}
where $\partial \Omega$ is the boundary of $\Omega$ and $L$ is a differential operator. Then the boundary conditions can be homogenized using
\begin{equation*}\label{bc0_r2}
u(\textbf{x})=v(\textbf{x})+h(\textbf{x}),
\end{equation*}
where $h$ satisfies the nonhomogeneous boundary conditions $\ Bu(\textbf{x})=g(\textbf{x})$.
After homogenization of the boundary conditions, the nonhomogeneous problem can be convert in the following form
\begin{equation*}\label{bc_r1}
 Lv(\textbf{x})=F(\textbf{x}),\hspace{.5 cm} \textbf{x}\in \Omega\subset R^{d},\ Bv(\textbf{x})=0,\hspace{.5 cm} \textbf{x}\in \partial \Omega
\end{equation*}
where $F(\textbf{x})=f(\textbf{x})-Lh(\textbf{x})$.
\begin{example}\label{ex1}
Consider the following linear two-point boundary value problem
\begin{equation}\label{e16}
  \left\{%
\begin{array}{ll}
    u''(x)=-\frac{\sinh(x)}{(1+\cosh(x))^{2}},\hskip .5 cm\hbox{$-1<x<1$} \\
    u(-1)=\alpha,u(1)=\gamma,
\end{array}%
\right.
\end{equation}
where $\alpha$ and $\gamma$ are given such that the exact solution
is $u(x)=\frac{\sinh(x)}{1+\cosh(x)}$.
\end{example}

\begin{figure}
\centering
\includegraphics[scale=.7]{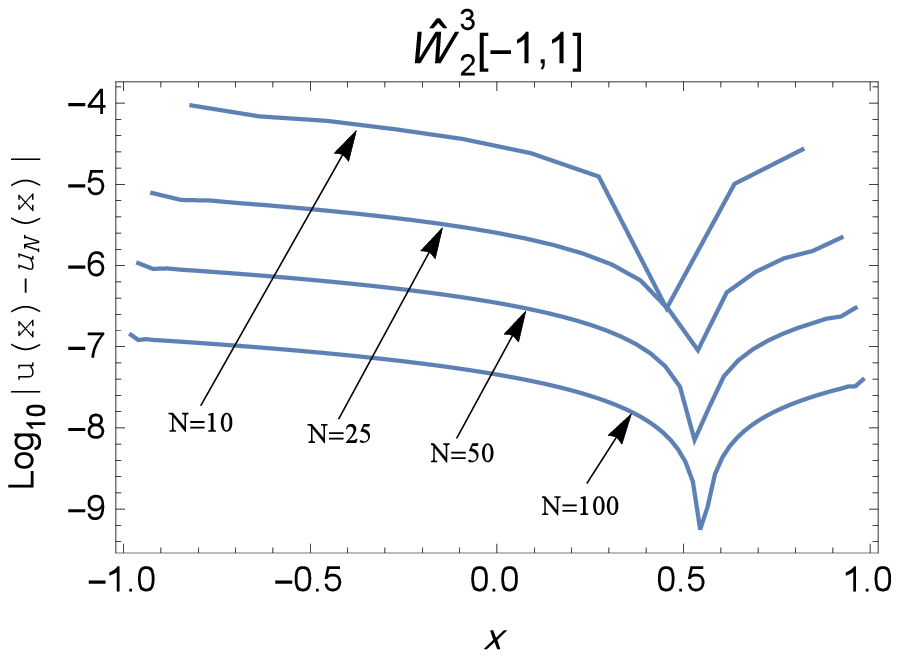}
\includegraphics[scale=.7]{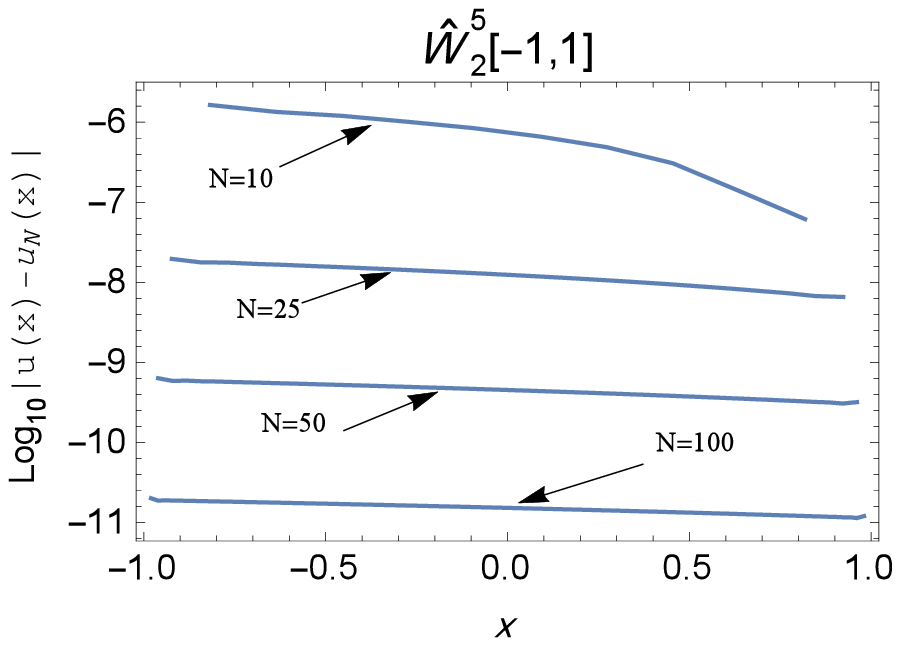}
\caption{\scriptsize{ Graph of $Log_{10}|u(x)-u_{N}(x)|$
with $m=3,5$ and $N=10,25,50,100$, for Example \ref{ex1}.}}\label{fig_r1}
\end{figure}

\begin{figure}
\centering
\includegraphics[scale=.7]{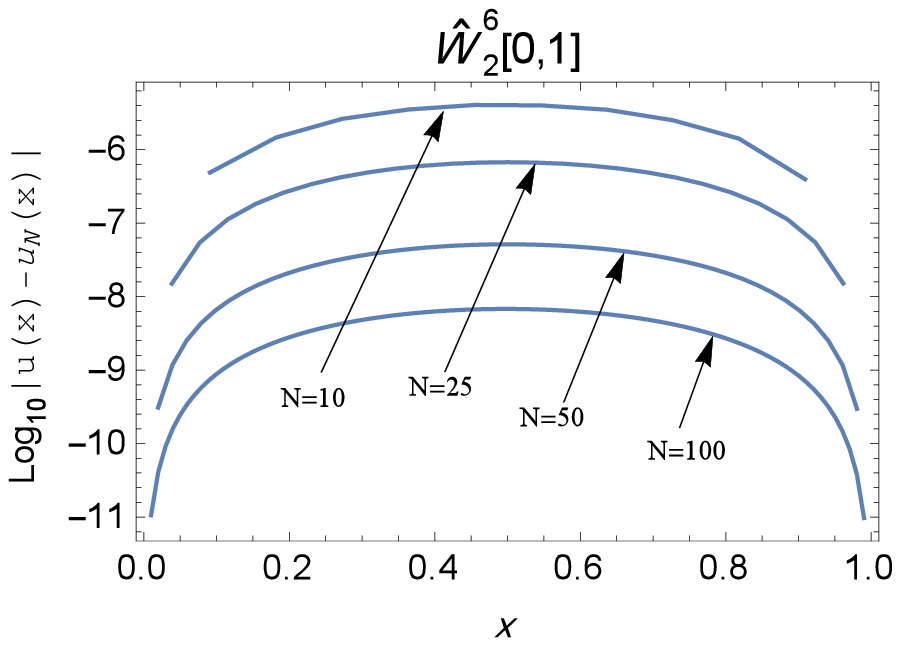}
\includegraphics[scale=.7]{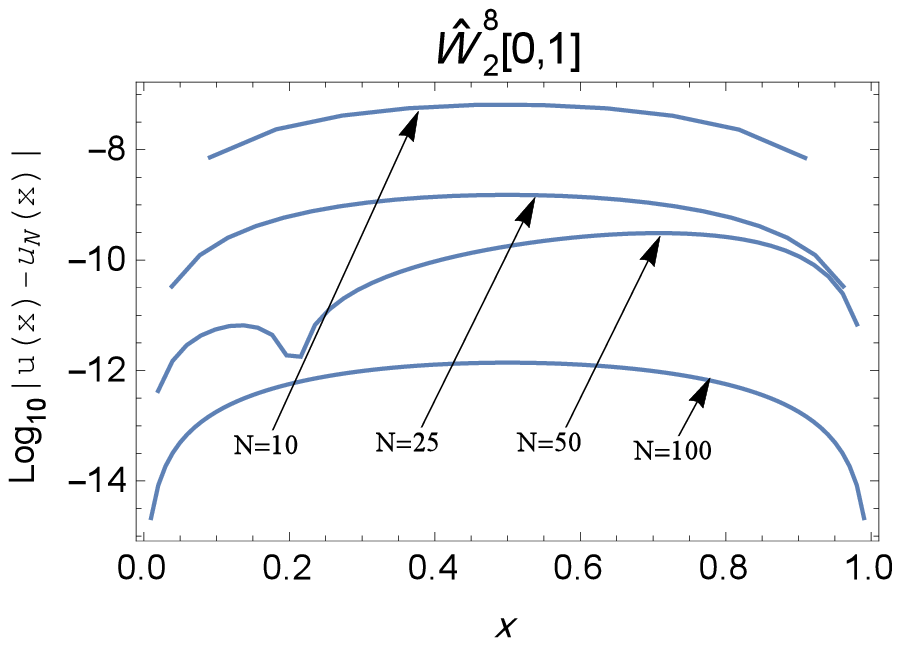}
\caption{\scriptsize{ Graph of $Log_{10}|u(x)-u_{N}(x)|$
with $m=6,8$ and $N=10,25,50,100$, for Example \ref{ex2}.}}\label{fig_r2}
\end{figure}

\begin{figure}
\centering
\includegraphics[scale=0.6]{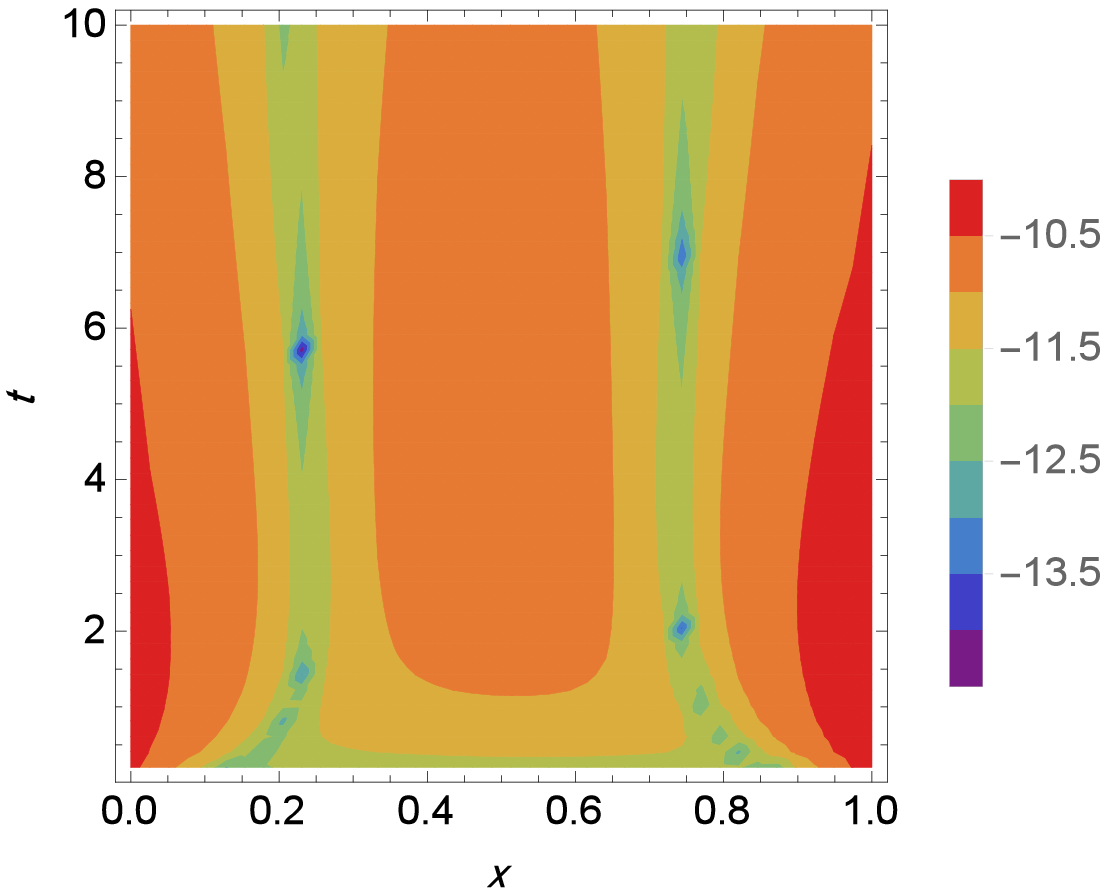}
\includegraphics[scale=0.6]{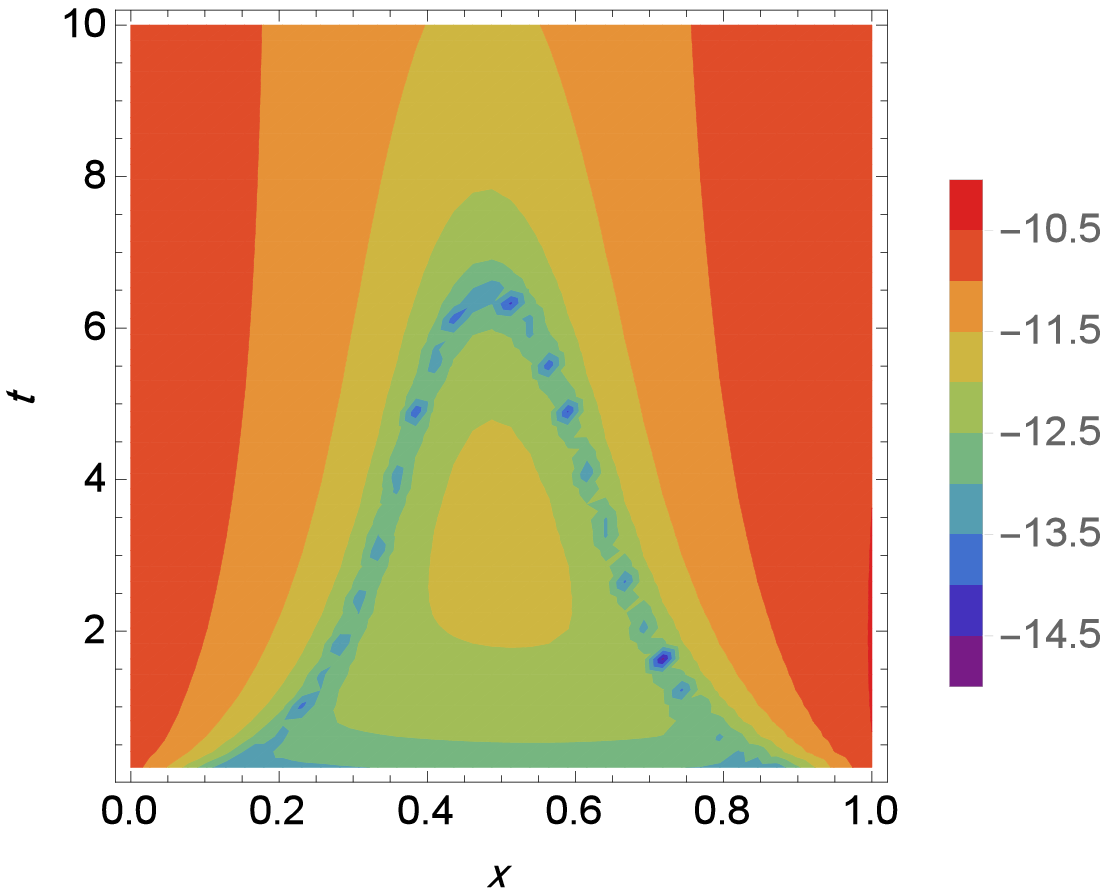}
\caption{\scriptsize{ Graph of $Log_{10}|u(x,t)-u_{N}(x,t)|$ in $(x,t)\in
[0,1]\times [0,10]$ with $v=0.01(left)$, $v=0.005(right)$, $\sigma=100,\Delta t=0.01$, $N=40$ and $\widehat{W}_{2}^{5}[0,1]$ for Example \ref{ex3}.}}\label{fig_r3}
\end{figure}

\begin{example}\label{ex2}
Consider the following fifth-order two-point boundary value problem
\begin{equation}\label{e016}
  \left\{%
\begin{array}{ll}
    u^{(5)}(x)+u(x)=g(x),\hskip .5 cm\hbox{$0<x<1$} \\
    u(0)=0,u(1)=0,\\
    u'(0)=1,u'(1)=-e,\\
    u^{(3)}(0)=-3,
\end{array}%
\right.
\end{equation}
where $g$ is given such that the exact solution is $u(x)=x(1-x)e^{x}$.
\end{example}
To solve Examples \ref{ex1} and \ref{ex2}, first we construct reproducing kernel
spaces $\widehat{W}_{2}^{m}[a,b]\subset W_{2}^{m}[a,b]$, where
$(m\geq 3)$ for Example \ref{ex1} and $(m\geq 5)$ for Example \ref{ex2} and in
which every function satisfies the homogenized boundary conditions. An approximate solution at the grid points $x_{i}$ might be
obtained by solving the discrete linear system
\begin{equation*}\label{e17}
L\mathbf{u}=\mathbf{f},
\end{equation*}
where $\mathbf{f}$ contains the values of function $f$ the righthand
function of differential equation after the homogenization, at the
grid points and $L$ is the differentiation matrix of differential
operator in \eqref{e16} and \eqref{e016}. In other words, the
solution at the grid points is given by
\begin{equation*}\label{e18}
\mathbf{u}=L^{-1}\mathbf{f},
\end{equation*}
and we see that invertibility of $L$ would be required. The maximum
absolute errors of approximate solutions of Example \ref{ex1} and
comparison with finite difference method and radial basis functions
collocation method are reported in Table \ref{tab_r1}. For radial
basis functions collocation method we used Gaussian kernel with
shape parameter $\epsilon=1$. The comparison of maximum absolute
errors, of Example \ref{ex2} with best reported results in
\cite{10,11,12,13} are shown in Tables \ref{tab001} and
\ref{tab0001}. In Figures~\ref{fig_r1} and \ref{fig_r2} we present
the maximum of absolute errors of approximate solutions in
logarithmic scale, for Example~\ref{ex1} and \ref{ex2} in different
reproducing kernel spaces and various values of N. The reported
results show that the accuracy of approximate solutions are closely
related to the smoothness order of the reproducing kernels and
values of N and as proved in \cite{babak3} more accurate approximate
solutions can be obtained using more mesh points and smoother
reproducing kernels. A Mathematica code of the implementation of the
method has been placed in the
http://www.abbasbandy.com/RKHS-COL.cdf.
\begin{example}\label{ex3}
Consider the Burgers' equation
\begin{equation}\label{e19}
  \left\{%
\begin{array}{ll}
    u_{t}+uu_{x}-vu_{xx}=0,\hskip .5 cm x\in(0,1),\ t\in(0,T], \\
   u(x,0)=f(x),\\
   u(0,t)=g_{1}(t),u(1,t)=g_{2}(t),
\end{array}%
\right.
\end{equation}
where $v=\frac{1}{Re}$ and $Re\geq 0$ is the Reynolds number
characterizing the size of viscosity and $f,g_{1},g_{2}$ is given
such that the exact solution is $u(x,t)=\frac{2v\pi
e^{-\pi^{2}vt}sin(\pi x)}{\sigma+e^{-\pi^{2}vt}cos(\pi x)}$, where
$\sigma$ is a parameter. For this example, the numerical results are
presented in Tables \ref{tab2} and \ref{tab3} for various values of
$N,v$ and they are compared with best reported results in
\cite{15,16,17} and finite difference method. Graphs of
$Log_{10}|u(x,t)-u_{N}(x,t)|$ in $(x,t)\in [0,1]\times [0,10]$ with
$v=0.01, 0.005$, $\sigma=100,\Delta t=0.01$, $N=40$ and
$\widehat{W}_{2}^{5}[0,1]$ are given in Figure \ref{fig_r3}. The
reported results show that more accurate approximate solutions can
be obtained using more mesh points and smoother reproducing kernels.
The numerical simulations show that the presented method is robust
and remain stable as time goes on. For solving time dependent
problems, we used differentiation matrix for the spatial
discretization together with an explicit Euler method with various
time steps as \cite{14}.
\end{example}

\begin{figure}
\centering
\includegraphics[scale=0.6]{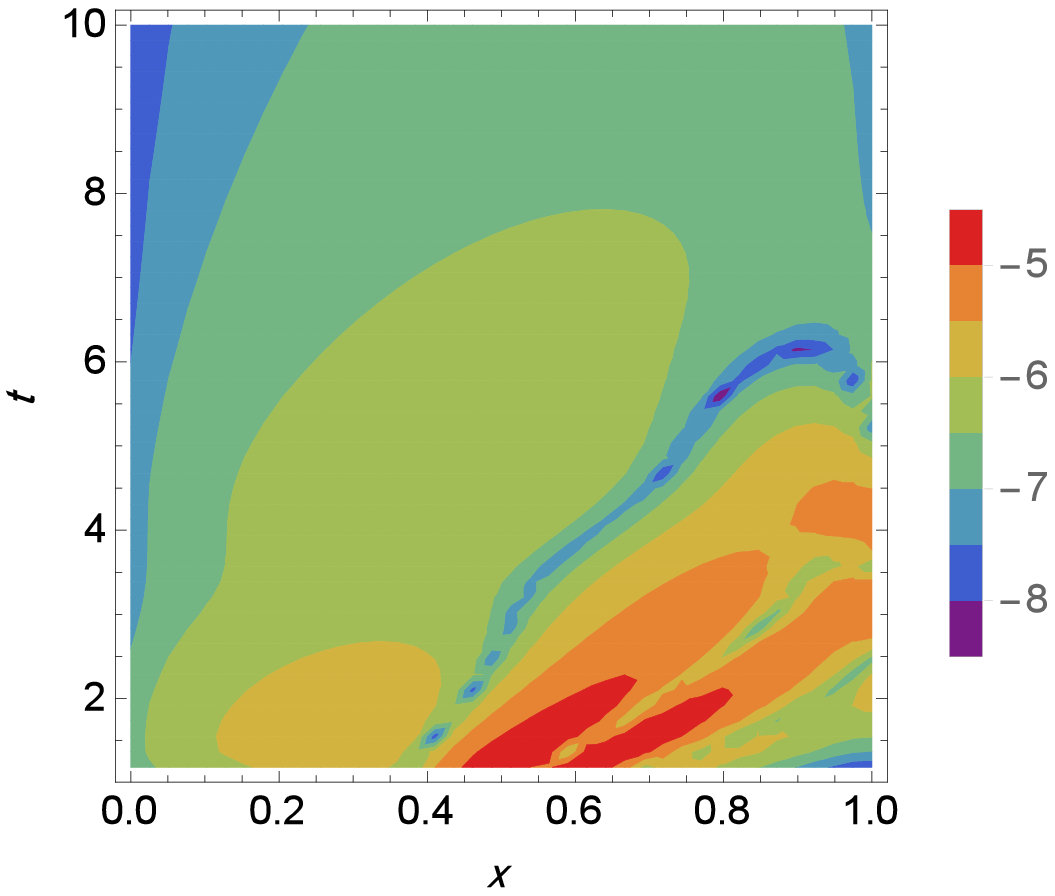}
\includegraphics[scale=0.6]{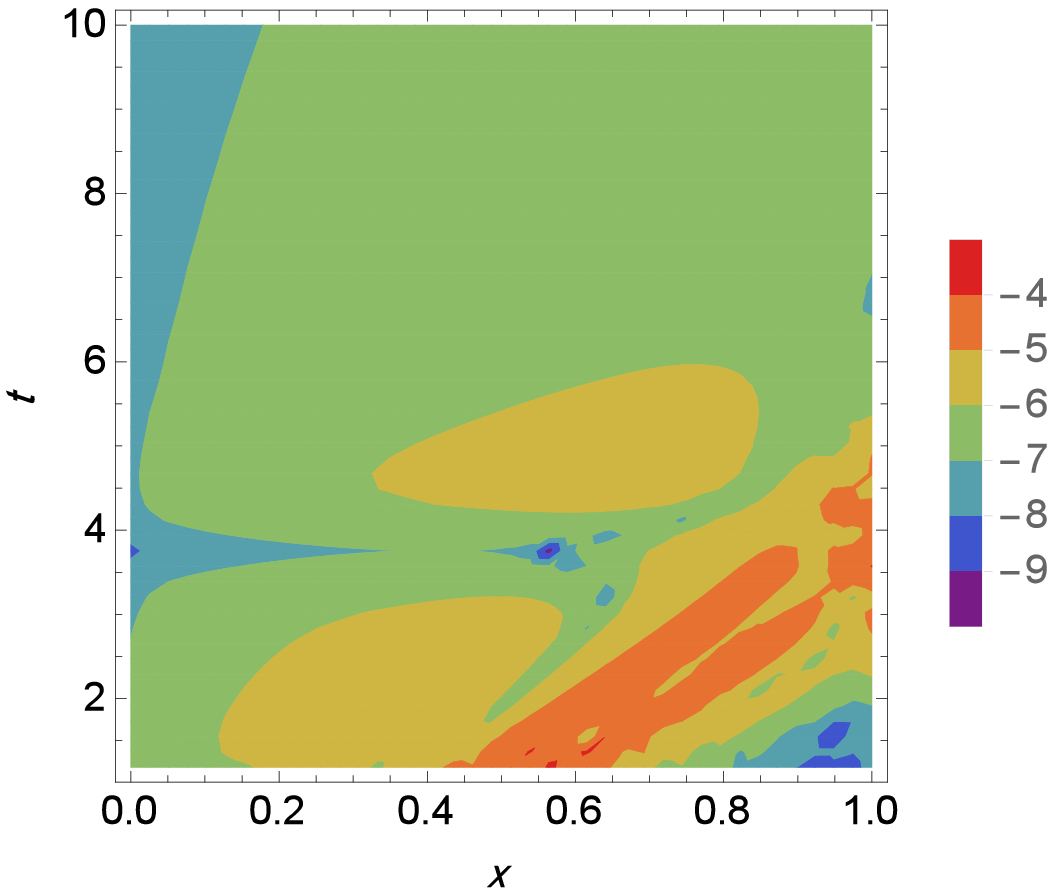}
\caption{\scriptsize{ Graph of $Log_{10}|u(x,t)-u_{N}(x,t)|$ in $(x,t)\in
[0,1]\times [1,10]$ with $v=0.01(left)$, $v=0.005(right)$, $\Delta t=0.01$, $N=40$ and $\widehat{W}_{2}^{5}[0,1]$ for Example \ref{ex4}.}}\label{fig_r4}
\end{figure}

\begin{example}\label{ex4}
Consider the Burgers' equation \eqref{e19} with $f,g_{1},g_{2}$ is
given such that the exact solution is
$u(x,t)=\frac{(\frac{x}{t})}{1+(\frac{t}{t_{0}})^{\frac{1}{2}}exp(\frac{x^{2}}{4vt})},t\geq
1$, where $t_{0}=exp(\frac{1}{8v})$. For this example, the numerical
results are presented in Table \ref{tab6} for various values of
$N,v$ and time $T$ and they are compared with best reported results
in \cite{18,19,20,21} and finite difference method. Graphs of
$Log_{10}|u(x,t)-u_{N}(x,t)|$ in $(x,t)\in [0,1]\times [1,10]$ with
$v=0.01, 0.005$, $\sigma=100,\Delta t=0.01$, $N=40$ and
$\widehat{W}_{2}^{5}[0,1]$ are given in Figure \ref{fig_r4}. We see
that the accuracy increases with increasing the mesh and smoothness
of kernels. The numerical simulations show that the presented method
is robust and remain stable as time goes on.
\end{example}

\begin{example}\label{ex5}
Consider the two dimensional Burger's equation,
\begin{equation}\label{e20}
u_{t}+uu_{x}+uu_{y}=vu_{xx}+vu_{yy},\hskip .5cm x_{0}\leq x \leq
x_{N},\hskip 0.1cm y_{0}\leq y \leq y_{N},\hskip 0.1cm t>0,
\end{equation}
with initial condition $u(x,y,0)=u_{0}(x,y)$ and the viscous
coefficient $v=\frac{1}{Re}>0$, $Re$ is the Reynolds number. The
Dirichlet boundary conditions is given such that the exact
solution is $u(x,y,t)=\frac{1}{1+e^{(x+y-t)/2v}}$, and
$x_{0}=y_{0}=0,x_{N}=y_{N}=1$.
\end{example}

\begin{theorem} \label{th002} \cite{22}
Let $W_{1}$ and $W_{2}$
be reproducing kernel spaces with reproducing kernels $K_{1}$ and
$K_{2}$. The direct product $\overline{W}=W_{1}\bigotimes W_{2}$ is
a reproducing kernel space and possesses the reproducing kernel
$\overline{K}(x_{1},x_{2},y_{1},y_{2})=K_{1}(x_{1},y_{1})K_{2}(x_{2},y_{2})$.
\end{theorem}

For solving multi-dimensional problems we are using the product of
reproducing kernels of $W_{2}^{m}[a,b]$ as kernels in
multi-dimensional domain. It is easy to see that these kernels are
strictly positive definite as proof of Proposition \ref{pr2}.

For Example \ref{ex5}, the comparison of Maximum absolute errors
$L_{\infty}$ and relative errors $L_{2}$ with radial basis functions
pseudospectral method and Chebyshev pseudospectral method using
$N=25$ and various $dt,m$ and $v$ in $\widehat{W}^{5}_{2}[0,1]$ are
presented in Table \ref{tab_r11}. For radial basis functions
pseudospectral method we used Gaussian kernel with optimal shape
parameter introduced in \cite{14}.

\begin{example}\label{ex6}
Consider the two dimensional Burger's equation
\begin{equation}\label{e21}
u_{t}+uu_{x}+uu_{y}=vu_{xx}+vu_{yy},\hskip .5cm x_{0}\leq x \leq
x_{N},\hskip 0.1cm y_{0}\leq y \leq y_{N},\hskip 0.1cm t>0,
\end{equation}
with initial condition $u(x,y,0)=u_{0}(x,y)$ and the viscous
coefficient $v=\frac{1}{Re}>0$, $Re$ is the Reynolds number. The
Dirichlet boundary conditions is given such that the exact solution
is $u(x,y,t)=0.5-tanh(\frac{x+y-t}{2v})$, and
$x_{0}=y_{0}=-0.5,x_{N}=y_{N}=0.5$. For Example \ref{ex6}, the
comparison of Maximum absolute errors $L_{\infty}$ and relative
errors $L_{2}$ with radial basis functions pseudospectral method and
Chebyshev pseudospectral method using $N=25$ and various $dt,m$ and
$v$ in $\widehat{W}^{5}_{2}[-0.5,0.5]$ are presented in Table
\ref{tab_r12}. For radial basis functions pseudospectral method we
used Gaussian kernel with optimal shape parameter introduced in
\cite{14}.
\end{example}

\begin{figure}
\centering
\includegraphics[scale=0.6]{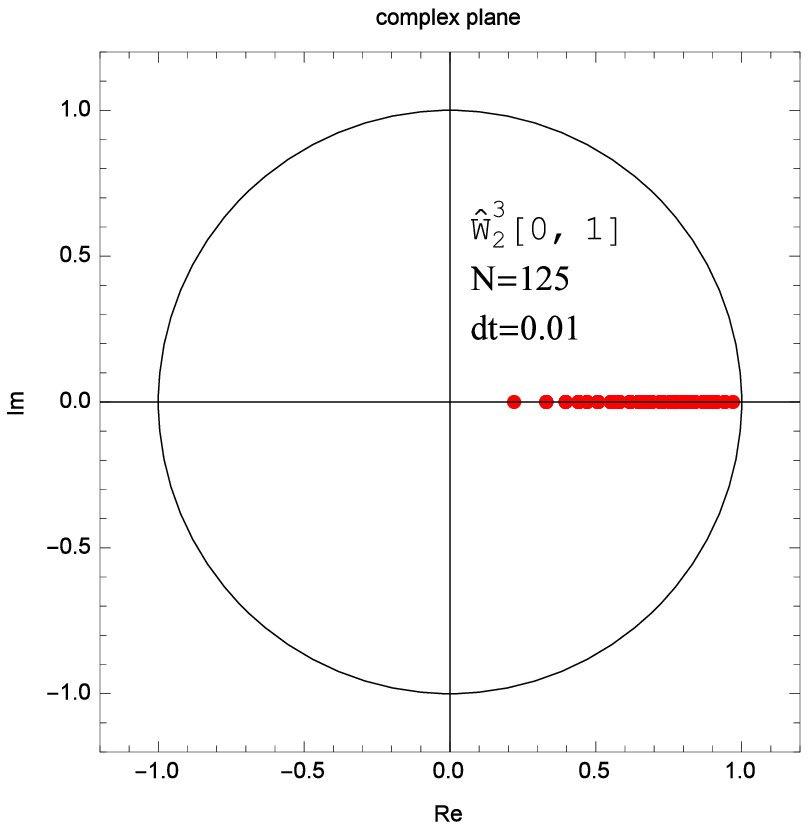}
\includegraphics[scale=0.6]{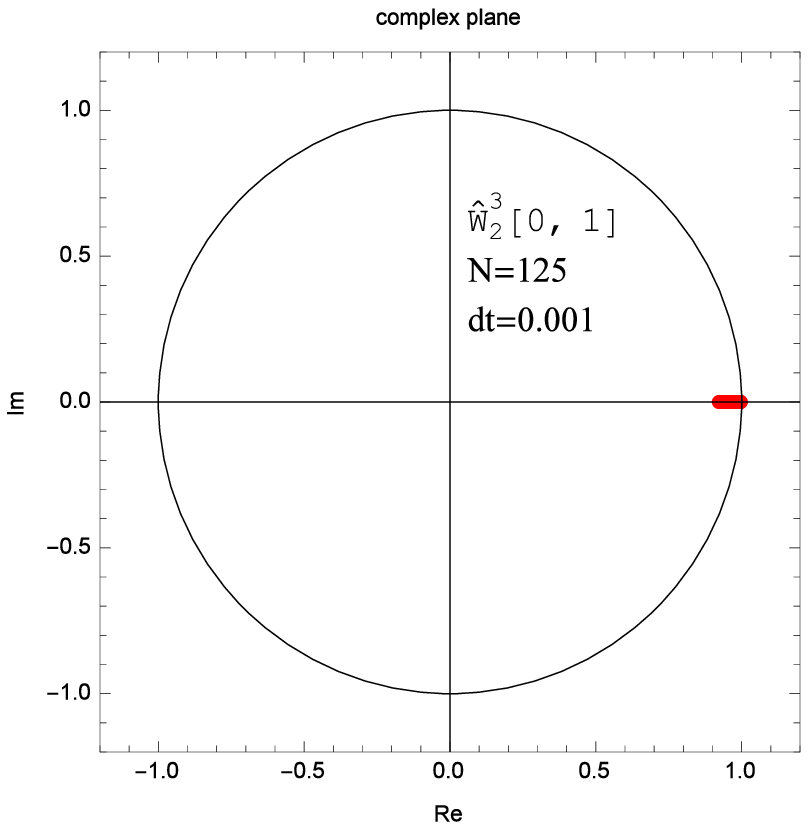}
\includegraphics[scale=0.6]{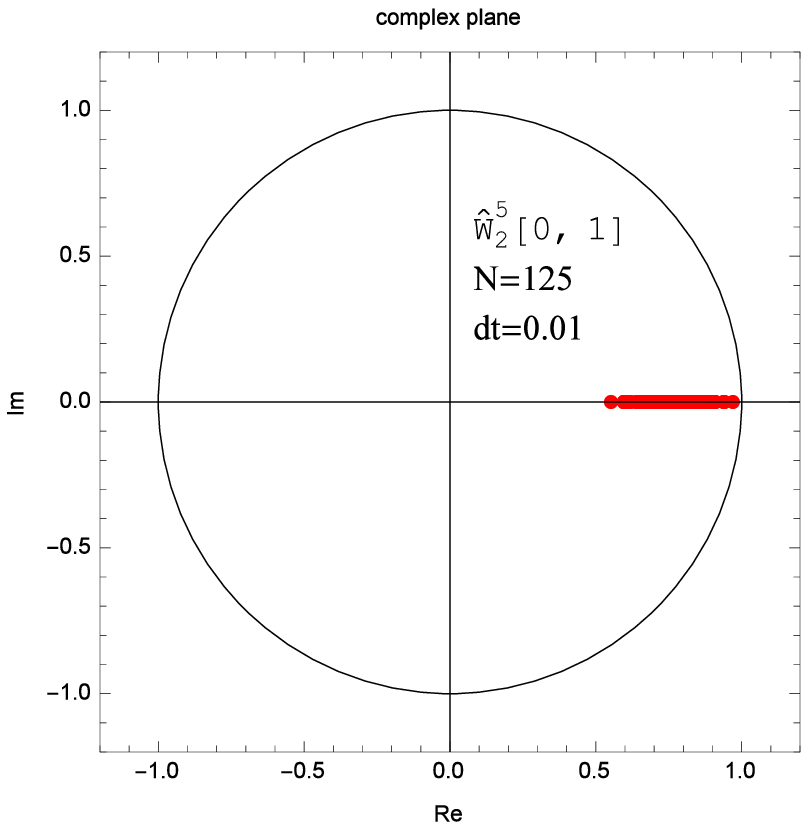}
\includegraphics[scale=0.6]{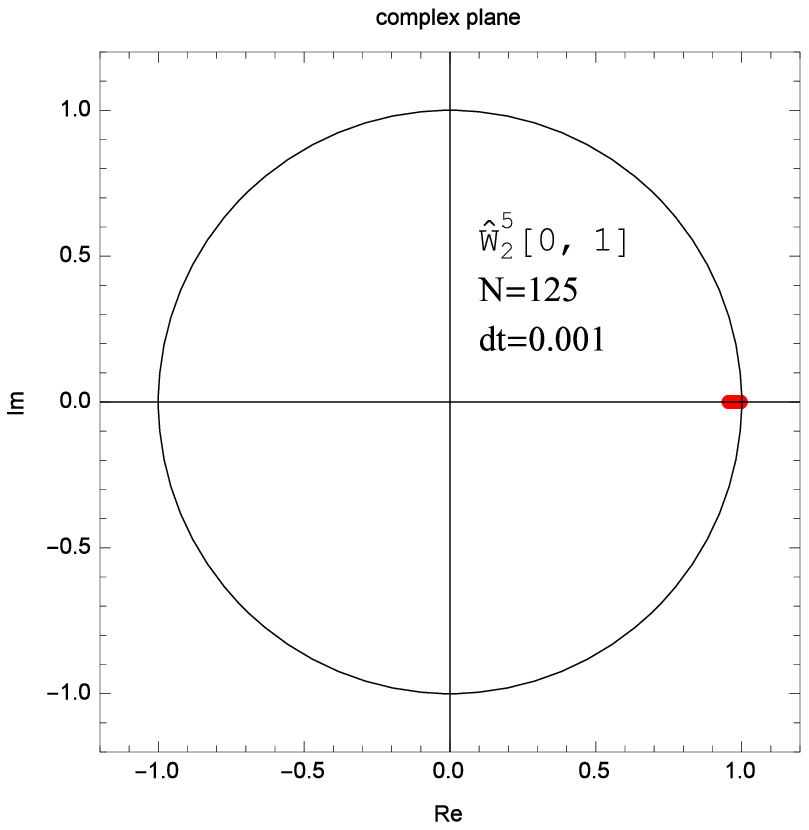}
\caption{\scriptsize{ Graph of eigenvalues of iteration matrix, with $N=125, \widehat{W}_{2}^{m}[0,1]$ and various $m,dt$ for Example~\ref{ex7}.}}\label{fig_r5}
\end{figure}

\begin{example}\label{ex7}
Consider the three dimensional problem
\begin{equation}\label{e22}
\frac{\partial u(x,y,z,t)}{\partial
t}=\frac{1}{\pi^{2}}\nabla^{2}u(x,y,z,t)-2 e^{t-\pi(x+y+z)},
\end{equation}
where $0\leq x,y,z\leq 1$ and $t>0$, with initial condition in
$t=0$, and Dirichlet boundary conditions which can be extracted
from the analytical solution,
$$u(x,y,z,t)=e^{t-\pi(x+y+z)}+x+y+z.$$
The relative errors $L_{2}$, of approximation solutions of Example
\ref{ex7} with various $T,dt,N$ and $m$ in
$\widehat{W}^{m}_{2}[0,1]$ are reported in Table~\ref{tab_r13}. In
Table~\ref{tab13} the numerical results are compared with best
reported results in \cite{23}. Figure~\ref{fig_r5} shows the graphs
of eigenvalues of iteration matrices of forward Euler time stepping,
in complex plane, with $N=125, \widehat{W}_{2}^{m}[0,1]$ and various
$m$ and $dt$ for Example~\ref{ex7}. The numerical simulations show
that the presented method is remain stable as time goes on, despite
of exponential growth of exact solution of problem in time and more
accurate approximate solutions can be obtained using more mesh
points and smoother reproducing kernels. From the results in Tables
\ref{tab_r13} and \ref{tab13}, we can see the implementation of
differentiation matrix with an explicit Euler method showed almost
identical behavior for a smaller time step so that we can be assured
that the inversion was indeed justified for this particular example.
\end{example}

\section{Conclusions}
In this paper, a new efficient meshless method, Reproducing kernel
Hilbert space mixed by meshless collocation method, based on
differentiation matrices which they are constructed by cardinal
functions of a RKHS, is proposed. In comparison with radial basis
function collocation method (Kansas method) we have the
nonsingularity in collocation matrix and Since the boundary
condition are imposed on trial space instead of collocation
matrix, the implementation of method is more simple and the method
is truly meshless. During the construction process we have
proposed a new and efficient algorithm to obtain the cardinal
functions of an RKHS and also we drive pointwise error estimate of
applying the differentiation matrices. To demonstrate the
computation efficiency, mentioned method is implemented for seven
examples and results have been compared with the reported results
in the literature which show the validity, accuracy and
applicability of the method.

\newpage

%

%
\clearpage

\begin{table}
 \centering
\scriptsize{
\begin{tabular}{ |c|c|c|c|c| }
  \hline
 N & N=10 & N=25 & N=50 &N=100\\
 \hline
 $\widehat{W}^{3}_{2}[-1,1]$&9.36088e-5  &  7.80543e-6  & 9.09414e-7 &1.40113e-7 \\
 \hline
 $\widehat{W}^{5}_{2}[-1,1]$&1.64341e-6&1.96221e-8  & 6.34336e-10&2.00868e-11 \\
  \hline
 Finite difference&8.14553e-5&1.3046e-5  &  3.27149e-6&8.17699e-7 \\
  \hline
 RBF collocation&6.32228e-5&4.85002e-6  &  9.44819e-7&5.60459e-7 \\
 \hline
\end{tabular}}
 \caption{\scriptsize{Maximum absolute errors of approximate solutions of Example \ref{ex1} and
 comparison with finite difference method and radial basis functions collocation method.}} \label{tab_r1}
\end{table}

\begin{table}
\centering
\scriptsize{\begin{tabular}{ |c|c|c|c| }
  \hline
 N & N=13 & N=26 & N=52 \\
 \hline
 \cite{11}, The fifth-order method & 1.3767e-4 & 7.1273e-6 &4.6950e-7\\
 \hline
\cite{11}, The seventh-order method &1.0024e-4 & 6.8397e-6& 4.4773e-7\\
 \hline
 \cite{10}&5.91739e-5& 3.40705e-7& 2.03387e-8\\
 \hline
Presented method, $\widehat{W}^{6}_{2}[0,1]$&4.14718e-6 & 3.29059e-7  & 4.60087e-8 \\
 \hline
 Presented method, $\widehat{W}^{8}_{2}[0,1]$&3.1921e-8 & 9.14844e-10  & 3.37252e-11 \\
 \hline
\end{tabular}}
\caption{\scriptsize{Maximum absolute errors, comparison
of results for Example \ref{ex2} with $N=13,26,52$.}}\label{tab001}
\end{table}

\begin{table}
\centering
\scriptsize{\begin{tabular}{ |c|c|c|c| }
  \hline
 N& N=10 & N=20 & N=40 \\
 \hline
 \cite{12} & 0.1570 & 0.0747 &0.0208\\
 \hline
\cite{13} &2.2593e-4 &1.3300e-5& 5.2812e-7\\
 \hline
 \cite{10}&6.29887e-5& 2.14116e-6& 7.00280e-8\\
 \hline
$\widehat{W}^{6}_{2}[0,1]$&4.06488e-6 & 6.75653e-7  & 9.77376e-8 \\
 \hline
$\widehat{W}^{8}_{2}[0,1]$&6.46414e-8 & 3.078e-9 & 1.19146e-10 \\
 \hline
\end{tabular}}
\caption{\scriptsize{Maximum absolute errors, comparison
of results for Example \ref{ex2} with $N=10,20,40$.}}\label{tab0001}
\end{table}


\begin{table}
\centering
\scriptsize{\begin{tabular}{ |c|c|c|c|c|c|c| }
  \hline
  N & \cite{15} &\cite{16} & \cite{17}& Finite difference& $\widehat{W}^{3}_{2}[0,1]$ & $\widehat{W}^{5}_{2}[0,1]$\\
  \hline
  10  & 1.2458e-7  &  1.215e-7  & 4.708e-8 & 3.2310e-6&2.18427e-7 &1.00476e-8 \\
  \hline
   20  &  3.3944e-8   & 3.062e-8   & 1.091e-8&1.5556e-6&5.0701e-8 &7.34209e-10 \\
  \hline
   40 & 1.1249e-8   & 7.644e-9   &  1.980e-9&7.6130e-7&8.45243e-9 &3.28094e-11\\
  \hline
\end{tabular}}
\caption{\scriptsize{Comparison of Maximum absolute
errors with existing numerical methods of Example \ref{ex3} for
$v=0.005,\sigma=100,\Delta t=0.01$ at $T=1.0$.}}\label{tab2}
\end{table}

\begin{table}
\centering
\scriptsize{\begin{tabular}{  |c|c|c|c|c|c|c| }
  \hline
  N &\cite{15} & \cite{16} & \cite{17}&Finite difference &  $\widehat{W}^{3}_{2}[0,1]$ & $\widehat{W}^{5}_{2}[0,1]$ \\
  \hline
  10  & 4.8808e-7  &  4.6280e-7  & 6.001e-11  &1.2250e-5&5.70664e-7&2.81404e-8 \\
  \hline
   20  &  1.4305e-7   & 1.1640e-7   & 1.010e-11&5.9232e-6&1.11397e-7&1.61939e-9 \\
  \hline
   40 & 5.6677e-8   & 2.9068e-8   &  1.277e-10&2.9100e-6&1.71283e-8 &6.57035e-11\\
  \hline
\end{tabular}}
\caption{\scriptsize{Comparison of Maximum absolute
errors with existing numerical methods of Example \ref{ex3} for
$v=0.01,\sigma=100,\Delta t=0.01$ at $T=1.0$.}}\label{tab3}
\end{table}


\begin{table}
\centering
 \scriptsize{\begin{tabular}{ |c|c|c|c|c|c|c|c| }
  \hline
  N&$\Delta t$&T &  \cite{18} & \cite{19} & Finite difference & $\widehat{W}^{3}_{2}[0,1]$&$\widehat{W}^{5}_{2}[0,1]$ \\
  \hline
  50&0.004&2.4  & 1.1e-3  &  4.0e-3  &7.8e-3& 3.11061e-5 &5.00091e-6  \\
  \hline
   100&0.001&2.4  &  2.8712e-4   & 9.9261e-4   &4.3e-3& 4.35295e-5&8.59652e-7 \\
  \hline
    N&$\Delta t$&T &  \cite{20} &  \cite{21}  &Finite difference &  $\widehat{W}^{3}_{2}[0,1]$&$\widehat{W}^{5}_{2}[0,1]$ \\
  \hline
  50&0.01&2.4  & 6.31491e-3  &  2.16784e-3  &8.5e-3& 4.24253e-5 &3.20613e-5   \\
  \hline
    50&0.01&1.8  & 5.12020e-3  &  2.47189e-3  &6.6e-3& 9.79958e-5 &6.82184e-5   \\
  \hline
\end{tabular}}
\caption{\scriptsize{Comparison of Maximum absolute
errors with existing numerical methods of Example \ref{ex4} for $v=0.005$.}}\label{tab6}
\end{table}

\begin{table}
\centering
 \scriptsize{\begin{tabular}{ |c|c|c|c|c|c| }
  \hline
  &  &&  RBF PS method & Chebyshev PS method&Presented method \\
  $dt$&T & v &  $L_{\infty}$\hskip 1cm$L_{2}$   &  $L_{\infty}$\hskip 1cm$L_{2}$& $L_{\infty}$\hskip 1cm$L_{2}$ \\
  \hline
  0.005 & 1 &  1& 9.0368e-3$\hskip .2cm$ 5.7021e-3  &  2.4851e-7$\hskip .2cm$ 1.9904e-7&4.25623e-9$\hskip .2cm$ 5.65924e-9 \\
  \hline
  0.001 & 1 &  1& 9.0366e-3$\hskip .2cm$ 5.7020e-3   & 9.6413e-8$\hskip .2cm$ 7.2229e-8& 4.09451e-9$\hskip .2cm$ 5.183e-9 \\
  \hline
  0.005&10&1& 8.195e-4$\hskip .2cm$ 2.9809e-4  & 4.8555e-7$\hskip .2cm$ 1.6113e-7&1.30473e-10$\hskip .2cm$ 8.08799e-11 \\
  \hline
  0.001&10&1& 8.1993e-4$\hskip .2cm$ 2.9822e-4  & 9.4330e-8$\hskip .2cm$ 3.2060e-8 &3.05613e-11$\hskip .2cm$ 2.10146e-11 \\
  \hline
  0.005 & 1 &  0.1& 2.9854e-1$\hskip .2cm$ 1.736e-1  & 4.27835e-3$\hskip .2cm$ 2.12004e-3&2.53845e-3$\hskip .2cm$ 2.572e-3 \\
  \hline
  0.001 & 1 &  0.1& 2.9816e-1$\hskip .2cm$ 1.735e-1   & 3.7823e-3$\hskip .2cm$ 2.0727e-3&2.83507e-3$\hskip .2cm$ 2.35878e-3 \\
  \hline
  0.005&10&0.1& 2.9024e-3$\hskip .2cm$ 1.1971e-3  & 3.6637e-15$\hskip .2cm$ 1.1098e-15&8.43362e-20$\hskip .2cm$ 3.299e-20 \\
  \hline
  0.001&10&0.1& 2.9024e-3$\hskip .2cm$ 1.1971e-3  & 2.0650e-14$\hskip .2cm$ 6.4789e-15&3.0511e-20$\hskip .2cm$ 1.58742e-20 \\
  \hline
\end{tabular}}
\caption{\scriptsize{Comparison of Maximum absolute
errors $L_{\infty}$ and relative errors
$L_{2}$ with existing numerical methods of Example~\ref{ex5} using $N=25$ and various $dt,m$ and $v$ in  $\widehat{W}^{5}_{2}[0,1]$.}}\label{tab_r11}
\end{table}

\begin{table}
\centering
 \scriptsize{\begin{tabular}{ |c|c|c|c|c|c| }
  \hline
  &  &&  RBF PS method & Chebyshev PS method&Presented method \\
  $dt$&T & v &  $L_{\infty}$\hskip 1cm$L_{2}$   &  $L_{\infty}$\hskip 1cm$L_{2}$& $L_{\infty}$\hskip 1cm$L_{2}$ \\
  \hline
  0.005 & 1 &  1& 5.4139e-2$\hskip .2cm$ 1.7775e-2  &  2.84257e-5$\hskip .2cm$ 1.01789e-5&6.68948e-6$\hskip .2cm$ 3.63904e-6 \\
  \hline
  0.001 & 1 &  1& 5.416e-2$\hskip .2cm$ 1.7781e-2  &  6.53996e-6$\hskip .2cm$ 2.41596e-6&1.62154e-6$\hskip .2cm$ 9.35968e-7 \\
  \hline
  0.005&10&1& 6.6292e-4$\hskip .2cm$ 1.8036e-4  &  1.43625e-8$\hskip .2cm$ 3.78170e-9&7.56945e-9$\hskip .2cm$ 3.26121e-9 \\
  \hline
  0.001&10&1& 6.6293e-4$\hskip .2cm$ 1.8036e-4  &  5.00597e-9$\hskip .2cm$ 9.10377e-10&1.52521e-9$\hskip .2cm$ 6.61975e-10 \\
  \hline
  0.005&1&0.1& 6.1243e-1$\hskip .2cm$ 1.5473e-1  &  1.25696e-2$\hskip .2cm$ 2.29108e-3&2.74444e-2$\hskip .2cm$ 7.36095e-3 \\
  \hline
  0.001&1&0.1& 6.1173e-1$\hskip .2cm$ 1.547e-1  &  1.23611e-2$\hskip .2cm$ 2.29390e-3&2.52867e-2$\hskip .2cm$ 6.71424e-3 \\
  \hline
  0.005&5&0.1& 6.58091e-3$\hskip .2cm$ 1.79631e-3  &  3.99680e-15$\hskip .2cm$ 6.59355e-16&6.66134e-16$\hskip .2cm$ 1.5666e-16 \\
  \hline
   0.001&5&0.1& 6.58245e-3$\hskip .2cm$ 1.79635e-3  &  1.22124e-14$\hskip .2cm$ 3.11285e-15&1.77636e-15$\hskip .2cm$ 4.13425e-16 \\
  \hline
\end{tabular}}
\caption{\scriptsize{Comparison of Maximum absolute
errors $L_{\infty}$ and relative errors
$L_{2}$ with existing numerical methods of Example~\ref{ex6} using $N=25$ and various $dt,m$ and $v$ in  $\widehat{W}^{5}_{2}[-0.5,0.5]$.}}\label{tab_r12}
\end{table}

\begin{table}
\centering
\scriptsize{\begin{tabular}{ |c|c|c|c|c|c| }
\hline
  T&dt & m &N=27 & N=64& N=125 \\
 \hline
 1& 0.01 & 3 &1.68375e-3&9.3954e-4&5.77729e-4\\
  \hline
 1& 0.01 & 5 &4.19018e-4&1.72518e-4&8.92225e-5\\
  \hline
 1 & 0.001 & 3  &1.65172e-3&9.1161e-4&5.5239e-4\\
 \hline
  1& 0.001 & 5  &3.93504e-4&1.41514e-4&6.18902e-5\\
 \hline
 5& 0.01 & 3  &3.10668e-2&1.51374e-2&8.47024e-3\\
 \hline
  5&0.01  & 5  &7.69788e-3&2.78871e-3&1.31209e-3\\
 \hline
 5& 0.001 & 3  &3.04834e-2&1.46877e-2&8.09662e-3\\
 \hline
  5&0.001  & 5  &7.22393e-3&2.28496e-3&9.06803e-4\\
 \hline
 \end{tabular}}
\caption{\scriptsize{Relative errors $L_{2}$, of
approximation solutions of Example \ref{ex7} with various $T,dt,N$ and $m$ in $\widehat{W}^{m}_{2}[0,1]$.}} \label{tab_r13}
\end{table}

\begin{table}
\centering
\scriptsize{\begin{tabular}{|c|c|c|c|c| }
 \hline
         &  $\widehat{W}^{3}_{2}[0,1]$,N=150 &  $\widehat{W}^{5}_{2}[0,1]$,N=150 &\cite{23},N=160  \\
 \hline
   dt    & $L_{\infty}$\hskip 1cm$L_{rms}$& $L_{\infty}$\hskip 1cm$L_{rms}$ &  $L_{\infty}$\hskip 1cm$L_{rms}$  \\
  \hline
  0.01  & 1.87124e-3\hskip .2cm 8.1899e-4 &2.56833e-4\hskip .2cm 1.27696e-4 &  2.98e-3\hskip .2cm 6.39e-4  \\
  \hline
   0.001 & 1.83218e-3\hskip .2cm 7.79576e-4 &2.09345e-4\hskip .2cm 8.43759e-5 & 3.88e-3\hskip .2cm 8.32e-4   \\
  \hline
 0.0001 & 1.82829e-3\hskip .2cm 7.75725e-4 &2.04977e-4\hskip .2cm 8.06016e-5 &2.84e-3\hskip .2cm 8.51e-4   \\
  \hline
\end{tabular}}
\caption{\scriptsize{$L_{\infty}$ and $L_{rms}$ , of approximation
solutions of Example \ref{ex7} at $T=1$ in $\widehat{W}^{3}_{2}[0,1]$ and $\widehat{W}^{5}_{2}[0,1]$.}} \label{tab13}
\end{table}

\clearpage

\end{document}